\documentclass[12pt]{article}
\usepackage{amsfonts,
amsmath,
latexsym,
amssymb,
makeidx,txfonts}
\usepackage[american]{babel}

\newcommand\RR{{\mathbb R}}

\newcommand\tr{{\rm tr\,}}
\newcommand\Tr{{\rm Tr\,}}
\newcommand\Det{{\rm Det\,}}

\newcommand\vol{{\rm vol\,}}

\newcommand\be{\begin{equation}}
\newcommand\ee{\end{equation}}
\newcommand\bea{\begin{eqnarray}}
\newcommand\eea{\end{eqnarray}}

\newtheorem{theorem}{Theorem}
\newtheorem{lemma}{Lemma}
\newtheorem{corollary}{Corollary}

\def\sideremark#1{\ifvmode\leavevmode\fi\vadjust{\vbox to0pt{\vss
\hbox to 0pt{\hskip\hsize\hskip1em
\vbox{\hsize2cm\tiny\raggedright\pretolerance10000
\noindent #1\hfill}\hss}\vbox to8pt{\vfil}\vss}}}

\begin{document}

\begin{titlepage}

\null
\vspace{3cm}

\centerline{\LARGE\bf Heat Kernel Asymptotics on Symmetric Spaces}
\bigskip
\bigskip
\bigskip
\centerline{\Large\bf\textbf{Ivan G. Avramidi}}
\bigskip
\centerline{\it Department of Mathematics}
\centerline{\it New Mexico Institute of Mining and Technology}
\centerline{\it Socorro, NM 87801, USA}
\bigskip

\centerline{\sc 
Dedicated to the Memory of Thomas P. Branson}
\medskip 

\begin{abstract}

We develop a new method for the calculation of the heat trace
asymptotics of the Laplacian on symmetric spaces that is based on a
representation of the heat semigroup in form of an average over the Lie
group of isometries and obtain a generating function for the whole
sequence of all heat invariants.

\end{abstract}



\end{titlepage}

\section{Introduction}
\setcounter{equation}0

The heat kernel is one of the most powerful tools in mathematical
physics and  geometric analysis (see, for example the books 
\cite{gilkey95,berline92,hurt83,avramidi00,kirsten01} and reviews 
\cite{avramidi87,camporesi90,avramidi96qftext,avramidi99,
avramidi02,vassilevich03}). In particular,  it is widely used to study
the propagators of quantum fields and the effective action in quantum
field theory  and the correlation functions and the partition function
in statistical physics. The short-time  asymptotic expansion of the
trace of the heat kernel determines the so-called spectral invariants of
the differential operator in question which  is intimately related to
the renormalization of quantum field theories
\cite{camporesi90,avramidi00,vassilevich03,kirsten01}, the
high-temperature expansion in statistical physics \cite{hurt83}, the
dynamics of integrable systems, in particular, the Korteweg-de Vries
hierarchy \cite{hurt83,avramidi00mn}, as well as spectral geometry and
index theorems \cite{gilkey95}. 

There has been a tremendous progress in the explicit calculation  of
spectral asymptotics in the last thirty years 
\cite{gilkey75,avramidi87,avramidi89,avramidi90,avramidi91,vandeven98,
yajima04}
(see also \cite{gilkey95,avramidi99,avramidi02,vassilevich03,kirsten01} 
and the references therein). However, due to the combinatorial explosion
in the complexity of the spectral invariants further progress in the
``brute force'' approach is  unlikely, even if employing computer
software for symbolic calculations. The results are so complicated that
it requires about twenty pages to describe them
\cite{vandeven98,yajima04}.
It seems that the further progress in the study of spectral asymptotics
can be achieved by restricting oneself to operators and manifolds with
high level of symmetry, for example, homogeneous and symmetric spaces,
which enables one to employ powerful algebraic methods. In some very
special  particular cases, such as group manifolds, spheres, rank-one
symmetric spaces, split-rank symmetric spaces etc,  it is possible to
determine the spectrum of the  differential operator exactly and to
obtain closed formulas for the  heat kernel in terms of the root vectors
and their multiplicities 
\cite{anderson90,camporesi90,dowker70,dowker71,hurt83,fegan83}. The
complexity of the method crucially depends on the global structure of
the symmetric space, most importantly its rank. Most of the results for
symmetric spaces are obtained for rank-one symmetric spaces only.
However, it is well known that the spectral asymptotics are  determined
essentially by local geometry. They are polynomial invariants in the
curvature with universal constants that do not depend on the global
properties of the manifold. It is {\it  this universal structure that we
are interested in this paper}. We will report on our results obtained in
the paper
\cite{avramidi96jmp} (see also \cite{avramidi94plb}).
Related problems in a more general context are discussed in
\cite{avramidi94grqc,avramidi95jmp36b,avramidi98wsp}.

\section{Spectral Asymptotics}
\setcounter{equation}0

\subsection{Laplacian}

Let $(M,g)$ be a smooth $n$-dimensional compact Riemannian manifold
without boundary  with a positive definite Riemannian metric $g$. Let
$TM$ and $T^*M$ be the tangent and cotangent bundles of the manifold
$M$. The metric $g$ defines in a natural way the Riemannian volume
element $d\vol$ on $M$. Let $C^\infty(M)$ be the space of smooth
real-valued functions on $M$. Using the invariant Riemannian volume
element  $d\vol$ on $M$ we define the natural $L^2$  inner product
in $C^\infty(M)$, and the Hilbert space $L^2(M)$ as the
completion of $C^\infty(M)$ in the corresponding norm.

The metric also defines the torsion-free compatible connection
$\nabla^{TM}$ on the tangent bundle $TM$,  so called Levi-Civita
connection. Using the Levi-Civita connection   we naturally obtain
connections on all bundles in the tensor algebra over $TM$ and $T^*M$;
the resulting connection will be denoted just by $\nabla$.  It
is usually clear which bundle's connection is being referred to, from
the nature of the section being acted upon. With our notation, Greek
indices, $\mu,\nu,\dots$, label the local coordinates $x=(x^\mu)$ on $M$
and range from 1 through $n$. Let $\partial_\mu$ be a
coordinate basis for the tangent space $T_xM$ at a point $x\in M$ and
$dx^\mu$ be dual basis for the cotangent space $T_x^*M$. We adopt the
notation that the Greek indices label the tensor components with respect
to local coordinate frame and range from 1 through $n$.
We also adopt the Einstein convention and sum over repeated indices.

Let $\nabla^*: C^\infty(T^*M)\to C^\infty(M)$ be the formal adjoint of 
$\nabla: C^\infty(M)\to C^\infty(T^*M)$ with respect to the $L^2$ inner
product. A partial differential operator $\Delta: C^\infty(M)\to
C^\infty(M)$ of the form
\be
\Delta=-\nabla^*\nabla=g^{\mu\nu}\nabla_\mu\nabla_\nu\,
\label{1ms}
\ee
is called the scalar Laplacian. In local coordinates the Laplacian is a
second-order partial differential operator of the form
\be
\Delta=|g|^{-1/2}\partial_\mu|g|^{1/2}g^{\mu\nu}\partial_\nu\,,
\ee
where $|g|=\det g_{\mu\nu}$. It is easy to see that
the principal symbol of the operator $(-\Delta)$ is 
\be
\sigma_L(-\Delta;x,\xi)=g^{\mu\nu}(x)\xi_\mu\xi_\nu\,,
\ee
where $\xi\in T^*_xM$ is a covector.

\subsection{Heat Kernel}

The following is known about the operator $(-\Delta)$
\cite{gilkey95}. First of all, it 
is elliptic and self-adjoint. More precisely, it is essentially
self-adjoint,  i.e. it has a unique self-adjoint extension to
$L^2(M)$. Second, the
operator $(-\Delta)$ has a positive definite leading symbol. The
spectrum of such operators forms a real nondecreasing sequence, 
$\{\lambda_k\}_{k=1}^\infty$, with each eigenspace being
finite-dimensional. The eigenvectors $\{\varphi_k\}_{k=1}^\infty$ are
smooth functions
and form an orthonormal basis in $L^2(M)$. Moreover, as $k\to
\infty$ the eigenvalues increase as $\lambda_k\sim Ck^2$ as
$k\to\infty$, with some positive constant $C$. These facts lead to the
existence of spectral asymptotics that will be discussed below.

The heat semi-group is a one-parameter family of bounded operators
$U(t)=\exp(t\Delta): L^2(M)\to L^2(M)$ for $t>0$. The integral kernel
$U(t|x,x')$ of this operator, called the heat kernel, is defined by
\be
U(t|x,x')
=\sum\limits_{k=1}^\infty e^{-t\lambda_k}\varphi_k(x)\varphi_k(x'),
\ee
where each eigenvalue is counted with multiplicities. It
satisfies the heat equation
\be
(\partial_t-\Delta)U(t|x,x')=0
\ee
with the initial condition
\be
U(0^+|x,x')=\delta(x,x')\,,
\ee
where $\delta(x,x')$ is the Dirac delta-function. {}For $t>0$ the heat
kernel is a smooth function on $M\times M$ with a well defined diagonal
\be
U^{\rm diag}(t|x)=U(t|x,x)\,.
\ee

\subsection{Spectral Invariants}

{}For any $t>0$ the heat semi-group $U(t)=\exp(t\Delta)$ is a
trace-class operator with a well defined $L^2$ trace, called the heat
trace,
\be
\Tr_{L^2}\exp(t\Delta)=\int\limits_M d\vol\; U^{\rm diag}(t)
=\sum\limits_{k=1}^\infty e^{-t\lambda_k}.
\ee

The heat trace is obviously a spectral invariant of the operator
$(-\Delta)$. It determines other spectral functions by integral
transforms. In particular,  the zeta-function, $\zeta(s,\lambda)$,  is
defined as the $L^2$ trace of the complex power of the operator
$(-\Delta-\lambda)$,
\be
\zeta(s,\lambda)=\Tr_{L^2}(-\Delta-\lambda)^{-s}
={1\over\Gamma(s)}\int\limits_0^\infty dt\; t^{s-1}\,
e^{t\lambda}\;\Tr_{L^2}\exp(t\Delta),
\ee
where $s$ and $\lambda$ are complex variables with  ${\rm
Re}\,\lambda<\lambda_1$ and ${\rm Re}\, s>n/2$.

The zeta function enables one to define, in particular, the regularized
determinant of the operator $(-\Delta-\lambda)$,
\be
\Det_{L^2}(-\Delta-\lambda)=\exp\left\{-
{\partial\over\partial s}\zeta(s,\lambda)\Big|_{s=0}\right\}\,,
\ee
which determines the one-loop effective action in quantum field theory.

\subsection{Asymptotic Expansion}

It is well known that the heat kernel diagonal  has the following
asymptotic expansion as $t\to 0^+$ \cite{gilkey95}
\be
U^{\rm diag}(t)\sim 
\sum\limits_{k=0}^\infty t^k a_k\,.
\label{600}
\ee
The coefficients $a_k$ are called local heat kernel coefficients. They
are scalar polynomials in the curvature and its covariant derivatives
which are  known explicitly  up to $a_5$. 
In particular,
\bea
a_0&=& 1\,,
\nonumber\\
a_1&=&{1\over 6}R\,,
\nonumber\\
a_2&=&\frac{1}{30}\Delta R+{1\over 72}R^2
-{1\over 180} R_{\mu\nu}R^{\mu\nu}
+{{1}\over {180}}R_{\mu\nu\alpha\beta}R^{\mu\nu\alpha\beta}\,,
\eea
where $R_{\mu\nu\alpha\beta}$ is the Riemann tensor,
$R_{\mu\nu}=R^\alpha{}_{\mu\alpha\nu}$ is the Ricci tensor and 
$R=R^{\mu\nu}{}_{\mu\nu}$ is the scalar curvature. The coefficient $a_3$
was computed in \cite{gilkey75}. The coefficient $a_4$ was first
computed in \cite{avramidi87} and published in
\cite{avramidi89,avramidi90,avramidi91} (see also \cite{avramidi00}).
The coefficient $a_5$ was computed in \cite{vandeven98,yajima04}.

The asymptotic expansion of the heat kernel diagonal can be
integrated over the manifold to give the asymptotic expansion
of the heat trace as $t\to 0$
\be
\Tr_{L^2}\exp(t\Delta)\sim 
(4\pi t)^{-n/2}\sum\limits_{k=0}^\infty t^k A_k\,,
\ee
where
\be
A_k=\int\limits_M d\vol\; a_k\,.
\ee
This is the famous Minakshisundaram-Pleijel asymptotic expansion. The
coefficients $A_k$ are spectral invriants of the Laplacian. They are
often called global heat kernel coefficients or
Hadamard-Minakshisundaram-De Witt-Seeley (HMDS) coefficients. This
expansion is of great importance in differential geometry, spectral
geometry, quantum field theory and other areas of mathematical physics,
such as the theory of Huygens' principle, heat kernel proofs of the
index theorems, Korteweg-De Vries hierarchy, Brownian motion etc. (see,
for example, \cite{hurt83}).

The general structure of the heat kernel coefficients can be described
as follows \cite{avramidi91,avramidi99,avramidi00,avramidi02}. 
We define symmetric tensors ${K}_{(j)}$ of type $(2,j)$ 
(that we call symmetric jets of order $(j-2)$) as
symmetrized covariant derivatives of the curvature
\be
{K}_{(j)}{}^{\alpha\beta}{}_{\mu_1\dots\mu_j}
=\nabla_{(\mu_1}\cdots\nabla_{\mu_{j-2}}
R^\alpha{}_{\mu_{j-1}}{}^\beta{}_{\mu_j)}\,,
\ee
where the parenthesis denote the complete symmetrization over all
indices included. 
Next, we define all possible scalar (orthogonal) invariants
of the form
\be
{J}^A_{m,{\bf n}}=\tr^A_g\; {K}_{(n_1+2)}\otimes\cdots\otimes
{K}_{(n_m+2)}\,,
\ee
where ${\bf n}=(n_1,\dots,n_m)$ is a multiindex, $|{\bf n}|
=n_1+\cdots+n_m$, and $\tr_g$ denotes contraction with the metric
to get a scalar. The index $A$ labels all possible contractions.
All invariants ${J}^A_{m,{\bf n}}$ have $m$ curvatures and
$|{\bf n}|$ derivatives of the curvatures. The classification
of these invariants is a separate interesting problem.

Then the local heat kernel coefficients have the grading
according to the number of derivatives,
that is, 
\be
a_k=\sum_{m=1}^k a_{k, m}\,,
\ee
where
\be
a_{k,m}
=\sum_{{{\bf n}\ge 0;\ |{\bf n}|=2k-2m}}\sum_{A} 
\,C^A_{\bf n}\,{J}^A{}_{m,{\bf n}}\,,
\ee
where $C^A_{\bf n}$ are some universal constants.
Notice that the dimension of the 
space of invariants ${J}^A{}_{m,{\bf n}}$
grows as a factorial for large $k$. This makes a
direct explicit calculation of the heat kernel coefficients
for large $k$ meaningless.

The leading terms in the heat kernel coefficients with the
highest number of derivatives were computed in
\cite{avramidi87,avramidi90b,avramidi91,branson90}.
They have the form
\be
A_{k}=\int\limits_M d\vol\;
\frac{(-1)^k k!}{2(2k+1)!}
\left\{2R^{\mu\nu}\Delta^{k-2}R_{\mu\nu}
+(k^2-k-1) R\Delta^{k-2}R\right\}
+\cdots\,,
\ee
where the dots denote the terms with less derivatives.
Notice that there are only two independent invariants.

We are interested in this paper in the opposite case, that is, the terms
\be
a_{k,k}=\sum_{A} \,C^A_{\bf 0}\,\tr^A_g\; 
{K}_{(2)}\otimes\cdots\otimes {K}_{(2)}\,,
\ee 
with no derivatives at all. More precisely, we assume that
the curvature is parallel, 
\be
\nabla_\mu R_{\alpha\beta\gamma\delta}=0\,,
\label{226}
\ee
in other words, we restrict ourselves to locally symmetric spaces. This
is a much more complicated case due to the presence of more invariants
and algebraic constraints on the curvature tensor. This is an
essentially non-perturbative calculation. Explicit results exist only in
some particular cases.

\section{Symmetric Spaces}
\setcounter{equation}0

We list below some well known facts from the theory of symmetric spaces
\cite{wolf72,ruse61,takeuchi91}. A Riemannian manifold with parallel curvature is called a 
locally symmetric space. A complete simply connected locally
symmetric space is called a globally symmetric space (or, simply, a 
symmetric space). A symmetric space is said to be of compact, noncompact
or Euclidean type if all sectional curvatures  are positive, negative or
zero. A product of symmetric spaces of  compact and noncompact
types is called a semisimple symmetric space. A general symmetric space is
a product of a Euclidean space and a  semisimple symmetric space.

It should be noted that our analysis is purely local.  We are looking
for a  universal local generating  function of the curvature invariants
in the category of locally symmetric spaces, that adequately reproduces
the  asymptotic expansion of the heat kernel diagonal.
This function should  give all the terms without covariant derivatives 
of the curvature $a_{k,k}$ in  the asymptotic expansion of the heat
kernel, in other words all heat  kernel coefficients $a_{k}$  for
any locally symmetric space. It turns  out to be much easier to obtain a
universal generating function whose Taylor coefficients
reproduce the heat kernel coefficients $a_{k}$ than  to compute them
directly. 

It is obvious that flat subspaces do not contribute to the heat
kernel coefficients $a_{k}$. Therefore, it is  sufficient to consider
only semisimple symmetric spaces. Moreover, since  the coefficients
$a_k$ are polynomial in the curvatures, one can restrict  oneself only
to symmetric spaces of compact type. Using the factorization  property
of the heat kernel and the duality  between the compact and  the
noncompact symmetric spaces one can obtain then the  results for the 
general case by analytic continuation. That is why we consider below
only compact symmetric spaces. 

\subsection{Holonomy, Isometry and Isotropy Groups}

Let $e_a=e^\mu_a\partial_\mu$ be a basis for the tangent space $T_xM$.
Let $e^a_\mu$ be the matrix inverse to $e^\mu_a$, defining the dual
basis $\omega^a=e^a_\mu dx^\mu\,$ in the cotangent space $T_x^*M$.  We
extend these bases by parallel transport along geodesics to get local
frames on the tangent bundle and the cotangent bundles, and, therefore,
on any tensor bundle, that are parallel along geodesics.  We adopt the
notation that the Latin indices from the beginning of the alphabet,
$a,b,c,\dots,$
label the tensor components with respect to this local frame and range
from 1 through $n$. Then the frame components of any parallel tensor
(such as the curvature tensor) are constant.
 
The components of the curvature tensor of a compact symmetric space
can be always presented in the form
\be
R_{abcd} = \beta_{ik}E^i{}_{ab}E^k{}_{cd}
\ee
where $E^i_{ab}$, $(i=1,\dots, p)$, with $p \le n(n-1)/2$ being a
positive integer, is a set of  
antisymmetric $n\times n$ matrices and $\beta_{ik}$
is some symmetric  nondegenerate positive definite 
$p\times p$ matrix. In the
following the Latin indices from the middle of the alphabet,
$i,j,k,\dots,$  will be
used to denote such matrices; they should be not confused with the Latin
indices from the beginning of the alphabet which denote tensor
components. The Latin indices from the middle of the alphabet will be
raised and lowered with the matrix $\beta_{ik}$ and its inverse
$(\beta^{ik})=(\beta_{ik})^{-1}$.

Next we define the traceless $n\times n$
matrices $D_i$, by $(D_i)^a{}_b=D^a{}_{ib}$,
where
\be
D^a{}_{ib}=-\beta_{ik}E^k{}_{cb}g^{ca}\,.
\ee
The matrices $D_i$ are known to be the generators of the  holonomy
algebra, ${\cal H}$, i.e. the Lie algebra of the restricted holonomy 
group, $H$,
\be
[D_i, D_k] = F^j{}_{ik} D_j\,,
\label{310}
\ee
where $F^j{}_{ik}$ are the structure constants of the holonomy group.
The matrix $\beta_{ik}$ plays the role of the metric of the holonomy group
with the scalar curvature 
\be
R_H=-\frac{1}{4}\beta^{ik}F^m{}_{il}F^l{}_{km}\,.
\ee
The structure constants define the traceless $p\times p$ matrices
$F_i$, by $(F_i)^j{}_k=F^j{}_{ik}$,
which generate the adjoint representation 
of the holonomy algebra,
\be
[F_i, F_k] = F^j{}_{ik} F_j\,.
\ee

Now, we go back to the equation (\ref{226}). It is an overdetermined
system of partial differential equations. By taking the commutator of
covariant derivatives we obtain the integrability condition of this
equation
\be
R_{fg e a}R^e{}_{b cd} 
-R_{fg e b}R^e{}_{a cd} 
+ R_{fg ec}R^e{}_{d ab}
-R_{fg ed}R^e{}_{c ab}
= 0\,.
\label{312}
\ee
By using the decomposition of the Riemann tensor introduced above 
we obtain from this equation
\be
E^i{}_{b c} D^c{}_{ka}
-E^i{}_{a c} D^c{}_{kb} 
= E^j{}_{a b} F^i{}_{j k}\,.
\label{313}
\ee
This is the most important equation that holds only in symmetric spaces;
it is this equation that makes a Riemannian manifold the  symmetric
space. 

To explore further consequences of the equation (\ref{313})
we introduce a new type of indices, the capital Latin indices,
$A,B,C,\dots,$ which split acording to
$A=(a,i)$ and run from $1$ to $N=p+n$.  We define a symmetric
nondegenerate  positive definite $N\times N$ matrix
\be
(\gamma_{AB}) = 
\left(
\begin{array}{cc}
g_{ab} & 0 \\
0 & \beta_{ik} \\
\end{array}
\right)\,.
\label{319}
\ee
This matrix and its inverse $(\gamma^{AB})=(\gamma_{AB})^{-1}$ will be used
to lower and to raise the capital Latin indices.
Next, we introduce a collection of new quantities $C^A{}_{BC}$ 
with the non-vanishing components
\be
C^i{}_{ab}=E^i{}_{ab}, \qquad 
C^a{}_{ib}=-C^a{}_{bi}=D^a{}_{ib}, \qquad 
C^i{}_{kl}=F^i{}_{kl}\,.
\label{317}
\ee
Let us also introduce rectangular $p\times n$ matrices 
$T_a$ by $(T_a)^j{}_c=E^j{}_{ac}$ and the $n\times p$ matrices
$\bar T_a$ by $(\bar T_a)^b{}_i=-D^b{}_{ia}$. 
Then we can define $N\times N$ matrices $C_A=(C_a,C_i)$
\be
C_a = \left(
\begin{array}{cc}
0 & \bar T_a \\
T_a & 0 \\
\end{array}
\right)\,,
\qquad
C_i = \left(
\begin{array}{cc}
D_i & 0 \\             
0 & F_i\\
\end{array}
\right),
\label{318}
\ee
so that $(C_A)^B{}_C=C^B{}_{AC}$. 

Then the matrices $C_A$ satisfy the commutation relations
\be
[C_A, C_B]=C^C{}_{AB}C_C\,,
\label{320}     
\ee
and generate the adjoint representation
of the Lie algebra ${\cal G}$ of some Lie group $G$
with the structure constants $C^A{}_{BC}$.
In more details, the commutation 
relations have the form
\bea
[C_a,C_b]&=&E^i{}_{a b}C_i,
\\{}
[C_i,C_a]&=&D^b{}_{ia}C_b,
\\{}
[C_i, C_k] &=& F^j{}_{i k} C_j\,,
\label{321}
\eea
which makes it clear 
that the holonomy algebra ${\cal H}$ is the subalgebra of the
Lie algebra ${\cal G}$. 

The matrix $\gamma_{AB}$ plays the role of the
metric on the group $G$ with the scalar curvature
\be
R_G=-\frac{1}{4}\gamma^{AB}C^C{}_{AD}C^D{}_{BC}\,.
\ee
It can be expressed in 
terms of the scalar curvature $R$ of the symmetric space $M$
and the scalar curvature $R_H$ of the isotropy subgroup $H $
\be
R_G = {3\over 4} R + R_H\,.
\label{429}
\ee

It is well
known that for compact symmetric spaces the group of isometries is
isomorphic to the Lie group $G$ defined in the previous section
(for more details see \cite{wolf72,avramidi96jmp}).  The
generators of isometries are the Killing vector fields
$(\xi_A)=(P_a,L_i)$, which form the Lie algebra of isometries
\be
[\xi_A,\xi_B]=C^C{}_{AB}\xi_C\,,
\label{343}
\ee
or
\bea
[P_a, P_b] &=& E^i{}_{a b} L_i\,,
\\{}
[L_i,P_a] &=& D^b{}_{ia} P_b\,,
\\{}
[L_i, L_k] &=& F^j{}_{i k} L_j\,.
\label{344}
\eea
The vector fields $L_i$ form the isotropy subalgebra of  the isometry
algebra, which is isomorphic to the holonomy algebra ${\cal H}$.
Thus, a compact symmetric space $M$ is isomorphic to the quotient
space of the isometry group by the isotropy subgroup $M=G/H$.

We will need the explicit form of the Killing vectors fields
in symmetric spaces. Let
us fix a point $x'$ in the manifold $M$. Let $d(x,x')$ be the geodesic
distance between a point $x$ and the fixed point $x'$ and 
\be
\sigma(x,x')=\frac{1}{2}[d(x,x')]^2\,.
\ee
Then the derivative $\nabla_\mu\sigma(x,x')$ 
is the tangent vector to the geodesic
connecting the points $x$ and $x'$ at the point $x$. We let
\be
y^a(x,x')=g^{ab}e^{\mu}_b(x)\nabla_{\mu}\sigma(x,x')\,
\ee
and
\be
K^a{}_b=R^a{}_{cbd}y^cy^d\,.
\ee
Notice that $y^a=0$ at $x=x'$.
Then one can choose the variables $y^a$ as new coordinates near $x'$
and show that
\cite{avramidi96jmp,avramidi87,avramidi91}
\be
P_a=\left(\sqrt{K}\cot\sqrt{K}\right)^b{}_a\frac{\partial}{\partial y^b}
\ee
\be
L_i=-D^b{}_{ia}y^a\frac{\partial}{\partial y^b}\,,
\ee
where $K$ is a $n\times n$ matrix with the entries $K^a{}_b$.


\section{Algebraic Methods for the Heat Kernel}
\setcounter{equation}0

\subsection{Heat Semigroup}

Let $k^A$ be the canonical coordinates on the isometry group $G$
so that each isometry is represented in the form
$\exp\left<k,\xi\right>$, where $\left<k,\xi\right>=k^A\xi_A$.
Then the left-invariant vector fields on the isometry group $G$ are
given by
\be
X_A =  X^M{}_{A} (k) {\partial\over\partial k^M}
\label{414}
\ee
where
\be
X^M{}_{A}(k) = \left({C(k)\over \exp C(k) -1}\right)^M{}_{A}\,.
\label{415}
\ee
and $C(k)=k^A C_A$.
The metric on the group $G$ is given by
\be
G_{MN} = \gamma_{AB} Y^{A}{}_{M} Y^{B}{}_{N}\,,
\label{418}
\ee
where $(Y^A{}_M)=(X^N{}_B)^{-1}$ is the inverse matrix to $X^N{}_B$.
Then it is easy to see that the determinant of the metric is
\be
|G|=\det G_{MN}
=|\gamma|\det{}_{{\cal G}}
\left({\sinh[C(k)/2]\over C(k)/2}\right)^2\,,
\label{419}
\ee
where $|\gamma|=\det \gamma_{AB}$.
Let $X_2$ be the Casimir operator on the group $G$ defined by
\be
X_2=\gamma^{AB}X_A X_B\,.
\ee

\begin{lemma}
Let $\Phi(t;k)$ be a function on the isometry group $G$
defined by
\be
\Phi(t;k)=(4\pi t)^{-N/2} 
\det{}_{\cal G}\left({\sinh[C(k)/2]\over C(k)/2}\right)^{1/2}
\exp\left\{ -{1\over 4t}\left<k,\gamma k\right>
+ {1\over 6} R_G t\right\}\,,
\label{46a}
\ee
where $\left<k,\gamma k\right>=\gamma_{AB}k^Ak^B$.
Then $\Phi(t;k)$ satisfies the heat equation
\be
\partial_t \Phi = |G|^{1/2}X_2|G|^{-1/2}\Phi\,,
\label{423}
\ee
and the initial condition
\be
\Phi(0;k)=|\gamma|^{-1/2}\delta(k)\,.
\label{428}
\ee

\end{lemma}

\noindent
{\bf Proof.}
First, we notice that the function $|G|^{-1/4}$ satisfies the
following equations
\be
X_2|G|^{-1/4} = {1\over 6} R_G |G|^{-1/4}\,,
\label{425}
\ee
and
\be
k^A{\partial\over \partial k^A} |G|^{-1/4} 
= {1\over 2} (N - X^A{}_A)
|G|^{-1/4}\,,
\label{426}
\ee
where
\be
X^A{}_A
={\rm tr}\,{}_{\cal G}\,{C(k)}\coth\left[{C(k)}\right]\,.
\label{427}
\ee
By using these equations we show that eqs. (\ref{423}) and (\ref{428}) hold
by a direct calculation.

One can show that the Laplacian on a symmetric space 
is simply the Casimir operator of the 
isometry group,
\be
\Delta = \gamma^{AB}\xi_A\xi_B\,,
\label{42}
\ee
and, therefore, belongs to the center of the enveloping algebra, i.e.,
\be
[\Delta, \xi_A] = 0\,.
\label{43}
\ee

\begin{theorem}
Let 
\be
\Psi(t) = \int\limits_{\RR^N} dk\; |\gamma|^{1/2}\Phi(t;k)
\exp\left<k,\xi\right>\,.
\label{49}
\ee
Then $\Psi(t)$ satisfies the heat equation
\be
\partial_t\Psi=\Delta\Psi
\label{411}
\ee
with initial condition
\be
\Psi(0)=1\,,
\label{412}
\ee
and, therefore,
\be
\Psi(t)=\exp(t\Delta)\,.
\label{417a}
\ee
This equation 
means that the formal power series as $t\to 0$ of both sides of
this equation are the same.
\end{theorem}

\noindent
{\bf Proof:}
We have 
\be
\partial_t\Psi(t)=\int\limits_{\RR^N} d k\;|\gamma|^{1/2} 
\partial_t\Phi(t;k)\exp\left<k,\xi\right>\,.
\label{422}
\ee
By using the previous Lemma we obtain
\be
\partial_t\Psi(t)=\int_{\RR^N} dk\; |\gamma|^{1/2}\exp\left<k,\xi\right>
|G|^{1/2}X_2|G|^{-1/2}\Phi(t;k)\,.
\label{421}
\ee
Now, by integrating by parts we get
\be
\partial_t\Psi(t)=\int_{\RR^N} dk\; |\gamma|^{1/2}\Phi(t;k)
X_2\exp\left<k,\xi\right>\,.
\label{421a}
\ee
Next, we show that
\be
X_B \exp\left<k,\xi\right>=\xi_B \exp\left<k,\xi\right>\,,
\label{413}
\ee
and, therefore,
\be
X_2 \exp\left<k,\xi\right>=\Delta\exp\left<k,\xi\right>\,.
\label{416}
\ee
Thus, the function $\Psi(t)$ satisfies the eq. (\ref{411}).
The initial condition (\ref{412}) for the function $\Psi(t)$
follows from the initial condition
(\ref{428}) for the function $\Phi(t;k)$.

\subsection{Heat Kernel Diagonal}

The heat kernel diagonal can be obtained by acting by the heat
semigroup on the delta-function,
\be
U^{\rm diag}(t;x)=\exp(t\Delta)\delta(x,x')\Big|_{x=x'}\,.
\ee
To be able to use integral representation for the heat semigroup
(\ref{49}) obtained above we need to compute the action of the
isometries $\exp\left<k,\xi\right>$ on the delta-function.

\begin{lemma}
Let $\omega^i$ be the canonical coordinates on the isotropy
group $H$ and $(k^A)=(q^a,\omega^i)$ be the natural splitting of the
canonical coordinates on the isometry group $G$.
Then
\be
\exp\langle k,\xi\rangle\delta(x,x')\Big|_{x=x'}
=|\eta|^{-1/2}
\det{}_{TM}\left(\frac{\sinh[D(\omega)/2]}{D(\omega)/2}\right)^{-1}
\delta(q)\,,
\label{424a}
\ee
where $D(\omega)=\omega^i D_i$ and $|\eta|=\det g_{ab}$.

\end{lemma}

\noindent
{\bf Proof.}
We choose the normal coordinates $y^a$ with the origin
at $x'$ defined above and consider
the equation of characteristics
\be
{d y^a\over ds} = \left(\sqrt {K(y)}\cot \sqrt
{K(y)}\right)^a_{\ b}q^b -\omega^iD^a{}_{ib}y^b\,.
\label{447}
\ee
By expanding the right hand side in the Taylor series we get
\be
{d y^a\over ds} = q^a
-\omega^iD^a{}_{ib}y^b
-\frac{1}{2}\beta^{ij}g_{gb}D^a{}_{ic}D^g{}_{jd} y^cy^dq^b
+O(y^3)\,.
\ee
Let $f^a(s,q,\omega)$ 
be the solution of the equation of characteristics with the initial condition
\be
f^a(0,q,\omega)=0\,.
\ee
Up to quadratic terms we obtain
\be
f^a(s,q,\omega)
=\left({1-\exp[-s D(\omega)]\over D(\omega)}\right)^a{}_{b}q^b
+O(q^2)
\label{452}
\ee
In particular, we find the Jacobian
\be
J(\omega)=\det\left({\partial f^a \over \partial q^b}
\right)_{q=0, s=1}
=\det{}_{TM}\left({\sinh[D(\omega)/2]\over D(\omega)/2}\right)\,.
\label{429a}
\ee
Then we have
\be
\exp\left<k,\xi\right>\delta(x,x')\Big|_{x=x'}
=|\eta|^{-1/2}\delta(f(1,q,\omega))\,.
\ee
By noticing that
\be
f^a(1,0,\omega)=0\,,
\label{442}
\ee
we finally obtain
\be
\delta(f(1,q,\omega)) 
=|\eta|^{-1/2}J(\omega)\delta(q)\,,
\label{441}
\ee
which proves the lemma.

One remark is in order here.
We implicitly assumed here that $q=0$ is the only
solution of the equation
\be
f^a(1,q,\omega)=0\,.
\ee
This is not necessarily true. This is the equation of closed
orbits of isometries and it has multiple solutions on
compact symmetric spaces. However, these global solutions,
which reflect the global topological structure of the manifold,
will not affect our local analysis. In particular, they do not
affect the asymptotics of the heat kernel. That is why, we have
neglected them here.

Now by using the above lemmas and the theorem we can compute the
heat kernel diagonal.

\begin{corollary}
The asymptotic expansion of the heat kernel diagonal as $t\to 0$
is given by the formal asymptotic expansion of the function
\bea
U^{\rm diag}(t)
&\sim &(4\pi t)^{-n/2}
\exp\left\{\left({1\over 8} R + {1\over 6} R_H \right)t\right\}
\int\limits_{\RR^p} \frac{d\omega}{(4\pi)^{p/2}}\;|\beta|^{1/2}
\exp\left\{-{1\over 4}\left<\omega,\beta\omega\right>\right\}
\nonumber\\
&&
\times\det{}_{\cal H}
\left({\sinh\left[\sqrt{t}\;F(\omega)/2\right]\over 
\sqrt{t}\;F(\omega)/2}\right)^{1/2}
\det{}_{TM}\left({\sinh\left[\sqrt{t}\;D(\omega)/2\right]\over 
\sqrt{t}\;D(\omega)/2}\right)^{-1/2}\,.
\label{437a}
\eea

\end{corollary}

\noindent
{\bf Proof.}
First, for the splitting $(k^A)=(q^a,\omega^i)$ we have $dk=dq\;d\omega$.
We compute the determinants of the metric (\ref{319})
\be
|\gamma|=|\beta|\;|\eta|\,,
\ee
where $|\beta|=\det\beta_{ik}$.
By using the equations (\ref{49}), (\ref{417a}),
and (\ref{424a}) and integrating over $q$ we obtain 
the heat kernel diagonal
\be
U^{\rm diag}(t)
=\int\limits_{\RR^p}d\omega\;|\beta|^{1/2}\Phi(t;0,\omega)
J(\omega)\,,
\ee
where $J(\omega)$ is given by (\ref{429a}).
Further, by using the eq. (\ref{318}) 
we compute the determinants
\be
\det{}_{\cal G}\left({\sinh[C(\omega)/2]\over C(\omega)/2}\right)
=\det{}_{TM}\left({\sinh[D(\omega)/2]\over D(\omega)/2}\right)
\det{}_{\cal H}\left({\sinh[F(\omega)/2]\over F(\omega)/2}\right)\,.
\label{437aa}
\ee
where $F(\omega)=\omega^i F_i$.
By using eqs. (\ref{46a}), (\ref{429a}) and (\ref{437a})
after scaling the integration variables
$\omega\to\sqrt{t}\;\omega$ we obtain finally
(\ref{437a}).

\subsection{Heat Kernel Asymptotics}

We introduce a Gaussian average over $\omega$ by
\be
\left<f(\omega)\right> = \int\limits_{\RR^p}
\frac{d\omega}{(4\pi)^{p/2}}\; |\beta|^{1/2}
\exp\left(-{1\over 4}\left<\omega,\beta\omega\right>\right)f(\omega)
\ee
Then 
\bea
U^{\rm diag}(t) &=&(4\pi t)^{-n/2}  
\exp\left\{\left({1\over 8} R + {1\over 6} R_H \right) t \right\}
\\
&&\times
\Bigg<
\det{}_{\cal H}\left({\sinh\left[\sqrt{t}\;F(\omega)/2\right]\over 
\sqrt{t}\;F(\omega)/2}\right)^{1/2}
\det{}_{TM}\left({\sinh\left[\sqrt{t}\;D(\omega)/2\right]\over 
\sqrt{t}\;D(\omega)/2}\right)^{-1/2}\Bigg>   
\nonumber
\eea

This equation can be used now to generate all heat kernel
coefficients $a_k$ for any locally symmetric space  simply by expanding
it in a power  series in $t$. By using the standard Gaussian averages
\bea
\left<\omega^i_1\cdots \omega^{i_{2k+1}}\right> &=& 0\,,
\label{440a}
\\
\left<\omega^{i_1}\cdots \omega^{i_{2k}}\right> 
&=& {(2k)!\over 2^{2k}k!}\beta^{(i_1 i_2}\cdots
\beta^{i_{2k-1}i_{2k})}
\label{441a}
\eea
one can obtain now all heat kernel coefficients in terms of various
contractions 
\be
\tr_\beta\left[F\otimes\cdots\otimes F
\tr_g\left(D\otimes\cdots\otimes D\right)\right]
\ee
of the matrices $D^a{}_{ib}$ and $F^{}j_{ik}$ with  the
matrices  $\beta^{ik}$ and $g^{ab}$. 
All these quantities are curvature invariants and
can be  expressed directly in terms of the Riemann tensor.

There is an alternative representation of the Gaussian average
in purely algebraic terms.
Let $b^j$ and $b^*_k$ be operators acting on a Hilbert space,
called creation and annihilation operators, 
that satisfy the following commutation relations
\be
[b^j,b^*_k]=\delta^j_k\,,
\ee
\be
[b^j,b^k]=[b^*_j,b^*_k]=0\,.
\ee
Let $|0\rangle$ be a unit vector in the Hilbert space, called the vacuum
vector, that satisfies the equations
\be
\langle 0|0\rangle=1\,,
\ee
\be
b^j|0\rangle=\langle 0|b_k^*=0\,.
\ee
Then the Gaussian average is nothing but the vacuum expectation
value 
\be
\langle f(\omega)\rangle
=\langle 0| f(b)\exp\langle b^*,\beta b^*\rangle |0\rangle\,,
\ee
where $\langle b^*,\beta b^*\rangle=\beta^{jk}b^*_jb^*_k$\,.
This should be computed by socalled normal ordering, that is, 
by simply commuting the operators $b_j$
through the operators $b^*_k$ until they hit the vacuum vector
giving zero. The remaining non-zero commutation terms precisely
reproduce the eqs. (\ref{440a}), (\ref{441a}).

\section*{Acknowledgements}

I would like to thank the organizers of the Midwest Geometry Conference
2006, in particular, Prof. Walter Wei, for the invitation and for the
financial support. I would like to dedicate this contribution to the
memory of Thomas P. Branson (1953-2006) whose sudden death shocked all
his friends and relatives. Tom was a good friend and a 
great mathematician who had a major
impact on many areas of modern mathematics.


\end{document}